\documentclass[11pt]{amsart}

\usepackage{amsmath,amssymb,latexsym}
\usepackage{amsfonts,amstext,amsthm,amscd}
\usepackage{graphicx}
\usepackage{amsthm}
\theoremstyle{plain}
\usepackage{color}
\usepackage[all]{xy}
\usepackage{mathtools}

\newtheorem*{tmaA}{Conjecture A}
\newtheorem*{conje}{Conjecture B}
\newtheorem{theorem}{Theorem}[section]
\newtheorem{cor}[theorem]{Corollary}

\newtheorem{ex}[theorem]{Example}
\newtheorem{lemma}[theorem]{Lemma}

\newtheorem{prop}[theorem]{Proposition}

\newtheorem{rest}[theorem]{Character Restrictions}
\numberwithin{equation}{section}
\def\bproof{\noindent{\textbf{Proof. }}}
\def\eproof{\noindent{\hfill $\blacksquare$}\bigskip}

\newcommand{\RR}{\mathbb{R}}
\newcommand{\CC}{\mathbb{C}}
\newcommand{\ER}{\widehat{\mathbb{C}}}

\newcommand{\NN}{\mathbb{N}}
\newcommand{\QQ}{\mathbb{Q}}

\newcommand{\pol}{\mathcal{P}ol(\CC)}
\newcommand{\rat}{\mathcal{R}at(\CC)}
\newcommand{\cng}{\mathcal{C}}

\begin{document}
\title{Elementary characters on semigroups: the rational case.}
\author{ Adri\'an Esparza-Amador \& Peter Makienko}
\begin{abstract}
Since polynomials form a subsemigroup of the semigroup of rational functions, every character on rational functions is a character on polynomials. On the other direction, not every character on polynomials is the restriction of a character on rational functions. What are the characters on polynomials that can be extended to rational functions? 

In this work, we conjecture that the only characters that can be extended are those that depends on the degree, often called \textit{elementary}. Also, we construct two example of character on polynomials, not elementaries, that cannot be extended to rational functions.
\end{abstract}
\maketitle

\section{Introduction}
A character over the semigroup of polynomials is an element of the set of homomorphisms from polynomials, with composition as binary operation, into the multiplicative group of complex numbers. The simplest non trivial example of a character over polynomials, is determined by the degree of a polynomial. Moreover, every multiplicative complex function of the degree function is again a character on polynomials. A question by E. Ghys considers the existence of non \textit{elementary} multiplicative character on polynomials, that is, a multiplicative complex function of degree function. In a recent paper, P. Makienko et. al. give an answer to such a question and non elementary characters over the affine group were constructed, which has an extension to the semigroup of polynomials. Since rational maps on the Riemann sphere are often considered as generalizations of polynomials, the notions of character can be extended to an algebraic study of rational functions as a semigroup, the semigroup of conformal endomorphism of the Riemann sphere. In this sense, every elementary character on polynomials can be extended to the semigroup of the rational functions. Is there any non elementary character over polynomial that can be extended to rational functions? 

Ritt's Theorem on prime decomposition of polynomials is a key on the construction of the non elementary character in \cite{cmp}. Although it is well know that there is no analogy of the Ritt's Theorem for rational functions, F. Pakovich in \cite{pak1} proves that the so called \textit{Laurent polynomials} is a semigroup for which the Ritt's Theorem holds. A couple of counterexamples to Ritt's Theorem (for rational functions) has been exhibit by Pakovich himself in \cite{pak} and by W. Bergweieler in \cite{ber}. One of such examples is used to construct a character on polynomials without extension to rational functions. Given the examples in this work, we consider the following conjecture. 
\begin{tmaA}
The only multiplicative character over polynomials that are restrictions of character on rational functions are elementaries.
\end{tmaA}
\textbf{Acknowledgements:} 
This project is supported by FORDECYT 265667 \emph{Program for a global and integrated advance of Mexican mathematics}.
\section{Preliminaries: semigroups}
Semigroups are among the most numerous objects in mathematics, and also among the most complex. Semigroup theory includes many subjects of interes: semigroups of transformation, semigroups associated to solutions of differential equations and in our interest case, the semigroup associated to rational complex iteration, to mention some examples. Below, we present some fundamentals in Semigroup Theory, for further information see \cite{gri}, \cite{how} and \cite{mag}.
\subsection{Definitions}
Given a set $S$ and a binary operation on $S$, $\cdot:S\times S\rightarrow S$, we say that the pair $(S,\cdot)$ is a \textbf{semigroup} if the binary operation is associative, that is, if
$$s\cdot(t\cdot u)=(s\cdot t)\cdot u,$$
for every $s,t,u\in S$. The binary operation is often denoted as a multiplication $x\cdot y=xy$.  We just say that $S$ is a semigroup when the binary operation is obvious.

Given a subset $T\subset S$ of a semigroup $S$, we say that $T$ is a \textbf{subsemigroup} of $S$ if $T$ is closed under the binary operation, that is, 
$$\text{if } t_1,t_2\in T, \text{ then }t_1\cdot t_2\in T.$$
An \textbf{identity element} of a semigroup $S$ is an element $e\in S$ such that 
$$ex=xe=x,\text{ for all }x\in S.$$
If an identity element exists for a semigroup, it is unique and usually denoted by $1$. If $S\neq\emptyset$ has no identity element, we may consider $S^1=S\cup \{1\}$, where $1\notin S$ and extend the binary operation on $S^1$ as follows: $1$ is an identity element and $xy$ is the same in $S$ and $S^1$ for all $x,y\in S$. If $S$ has an identity element then $S=S^1$. 

Similarly, a \textbf{zero element} of a semigroup $S$ is an element $z\in S$ such that
$$S\neq\{z\}\text{ and }zx=xz=z\text{ for all }x\in S.$$
When a zero element exists for a semigroup $S$ it is unique and usually denoted by $0$. If a non empty semigroup $S$ has no zero element, consider $S^0=S\cup\{0\}$, where $0\notin S$, and extend the binary operation on $S^0$ as follows: $0$ is a zero element and $xy$ is the same in $S$ and $S^1$ for all $x,y\in S$. When $S$ has a zero element then $S=S^0$. 

An element $e\in S$ is called \textbf{idempotent} if $ee=e^2=e$. If $A$ and $B$ are subsets of a semigroup $S$, then by definition $AB$ is the subset
$$AB=\{ab\in S|a\in A,\ b\in B\}$$
of $S$. In this setting, a subset $I\subset S$ is called a \textbf{left ideal} (respectively \textbf{right ideal}) of $S$ if $SI\subset I$ (resp. $IS\subset I$). $I$ is an \textbf{ideal} if is both left and right ideal. A mapping $\phi:S\rightarrow T$, between semigroups $S$ and $T$ is called a \textbf{homomorphism} if it preserves products: 
$$\phi(xy)=\phi(x)\phi(y) \text{ for all } x,y\in S, $$
where the product at the left of above equality is on $S$ while the product at the right side is on $T$.

An equivalence relation $\mathcal{C}$ on $S$ that admits multiplication ($a\ \cng\ b\ \&\ c\ \cng\ d \ \Rightarrow ac\ \cng\  bd$) is called a \textbf{congruence}. The resulting semigroup $S/\cng$ is the \textbf{quotient} of $S$ by $\cng$. 
\begin{prop}\label{rees}
When $I$ is an ideal of a semigroup $S$, the relation $J$ defined by
$$a\ J\ b \Leftrightarrow \ a=b \text{ or }a,b\in I$$
is a congruence on $S$. 
\end{prop}
The congruence $J$ in Proposition \ref{rees}, is called the \textbf{Rees congruence} of the ideal $I$ and the quotient semigroup $S/J=S/I$ the \textbf{Rees quotient} or the \textbf{Rees factor} of $S$ by $I$. Note that the Rees quotient can be expressed as 
$$S/J=S/I=\{I\}\cup\{x\in S|x\notin I\}.$$
Then, $S/I$ can be regarded as consisting of $I$ (as an element) together with the elements not in $I$ ($S\backslash I$). If $I\neq S,\emptyset$, the $J$-class $I\in S/I$ may be denoted by $0$ and then $S/I=(S\backslash I)\cup\{0\}$, defining the binary operation 
$$x*y=xy\text{ if }x,y\notin I \text{ and } x*y=0\text{ otherwise},$$
the quotient $Q=S/I$ is a semigroup with zero. 
\subsection{Multiplicative characters}
Given a semigroup (or a group) $G$ and a field $F$, a \textbf{character} of $G$ in $F$ is an homomorphism $\varphi:G\rightarrow F$. We say that the character is \textbf{multiplicative} or \textbf{additive} if the homomorphism on $F$ is with respect to the product or the addition on $F$ respectively. 
\begin{ex}\label{ex1}
Take $G=Aff(\CC)$, the group (under composition) of conformal automorphism of the complex plane, with $F=\CC$, and define $$\varphi(az+b)=a.$$ Since $(az+b)\circ (cz+d)=acz+ad+b$, we have $$\varphi((az+b)\circ (cz+d))=a\cdot c=\varphi(az+b)\cdot\varphi(cz+d),$$ this way, $\varphi$ is a multiplicative character. 
\end{ex}
We are interested in the case that $G$ is a semigroup, more specifically the semigroup of complex polynomials $Pol(\CC)$ with composition as semigroup operation and $F$ is the multiplicative group $\CC$. That is, $\varphi:Pol(\CC)\rightarrow \CC$ and satisfies
$$\varphi(P_1\circ P_2)=\varphi(P_1)\cdot \varphi(P_2),$$
for all $P_1,P_2\in Pol(\CC)$. 
The trivial examples of multiplicative characters in $Pol(\CC)$ are the constant maps $\varphi(P)=c$ for every $P\in \pol$, where $c\in\{0,1\}$. The basic non trivial example of a multiplicative character in $\pol$ is the degree function,
$$\varphi(P)=\text{deg}(P),$$
for all $P\in\pol$. In fact, every multiplicative complex function $f:\CC\rightarrow\CC$ generates a new multiplicative character in $\pol$: $\varphi_f(P)=f(\text{deg}(P))$. This type of characters, multiplicative complex functions of the degree function, are called \textbf{elementary}. 

In \cite{cmp}, authors construct a multiplicative character on $\pol$ which is not elementary. This way, they give an answer to the question of E. Ghys. Is there a multiplicative character on $\pol$ which is not elementary? The construction is based on the \textbf{Decomposition Ritt's Theorem} for polynomials, see \cite{ritt1}. 

Given a polynomial $P\in\pol$, we say that $P$ is \textbf{prime} or \textbf{indecomposable} if always that $P=Q\circ R$, then $\text{deg}(Q)=1$ or $\text{deg}(R)=1$. Otherwise, $P$ is called decomposable. If $P$ is decomposable, we say that a decomposition $P=P_1\circ...\circ P_n$ is a \textbf{prime decomposition} if each $P_i$ is a prime polynomial. We can state now the first part of Ritt's Theorem, \cite{ritt1}. 
\begin{theorem}
Let $P=P_1\circ...\circ P_n$ and $P=Q_1\circ...\circ Q_m$ be two prime decompositions of a polynomial $P$, then $m=n$. 
\end{theorem}
Then, the length of a prime decomposition of a given polynomial is a well define map $l:\pol\rightarrow \CC$. Moreover, again by the Ritt's Theorem, $l$ is an additive homomorphism, and then an additive character. Hence, defining the map $e^l:\pol\rightarrow\CC$ given by 
$$e^l(P)=\exp(l(P)),$$
then, $e^l$ is a multiplicative character which do not depends on the degree of the polynomial, and then it is not an elementary character. 

In general, every additive character can be turned up into a multiplicative character via the exponential map.

Also, using the second part of Ritt's Theorem, authors in \cite{cmp} give a method to construct characters in $\pol$. 
\begin{theorem}\label{defcha}
Let $\phi$ be a complex function, defined on the set of prime polynomials satisfying:
\begin{itemize}
\item[(i)] $\phi(c)=0$ for every constant $c$;
\item[(ii)] if $P_1,P_2, P_3, P_4$ are prime polynomials with $P_1\circ P_2=P_3\circ P_4$, then 
$$\phi(P_1)\cdot\phi(P_2)=\phi(P_3)\cdot\phi(P_4).$$
\end{itemize}
Then $\phi$ generates a multiplicative character $\Phi$. Conversely, if $\Phi$ is a multiplicative character in $\pol$, which is not the constant map $1$, then $\Phi$ satisfies (i) and (ii) above. 
\end{theorem} 
Given a polynomial $P$ with prime decomposition $P=P_1\circ...\circ P_n$, then $\Phi$ is defined as 
$$\Phi(P)=\Phi(P_1\circ...\circ P_n)=\phi(P_1)\cdot...\cdot \phi(P_n).$$
\subsection{Rational functions}
Since the set of rational functions $R:\ER\rightarrow\ER$ is also a semigroup with composition as semigroup operation, we can think of multiplicative characters of the semigroup of complex rational functions, $\rat$. 

Again, we have trivial multiplicative characters given by the constant maps $\varphi:\rat\rightarrow c$, $c\in\{0,1\}$, and the non trivial multiplicative character given by the degree function $\varphi(R)=\text{deg}(R)$, $R\in\rat$. Also, as in the case of polynomials, the degree function generates multiplicative characters through multiplicative complex functions, called again \textbf{elementary characters}. 

As in the example in \cite{cmp} of Affine characters for polynomials, it is possible to construct non constant characters as extensions through the group of M\"obius transformations $PSL(2,\CC)$. 
\begin{ex}[Fractional Linear Characters]
Let $H$ be the ideal of non injective rational functions. Any multiplicative character $\varphi:PSL(2,\CC)\rightarrow \CC$ admits an extension to a multiplicative character in $\rat$. For instance, put $\varphi(c)=0$ for every constant $c$, and $\varphi(h)=0$ for $h\in H$.
\end{ex}
This type of extensions (defining $\varphi\equiv0$ outside a subsemigroup or an ideal of $\rat$) are often called trivial extensions. 

Now, since $\pol\leq\rat$, the semigroup of complex polynomials is a subsemigroup of the complex rational functions, every restriction of a multiplicative character in $\rat$ to $\pol$ is a multiplicative character in $\pol$, that is, if $\varphi:\rat\rightarrow\CC$ is a multiplicative character in $\rat$, then $\varphi\large |_{\pol}:\pol\rightarrow\CC$ is a multiplicative character in $\pol$. 

A natural question arise: given a non trivial multiplicative character in $\pol$, $\varphi:\pol\rightarrow\CC$, can $\varphi$ be extended to a multiplicative character in $\rat$? That is, there exist a multiplicative character $\Phi:\rat\rightarrow\CC$ such that $\Phi\large|_{\pol}=\varphi$? Note that elementary characters on polynomials naturally extend to rational functions. 
\subsection{Characters and graduations} 
We have seen that an equivalence relation on a semigroup $S$ compatible with the product in $S$ is called a congruences. A canonical example of a congruence in the semigroup of polynomials, is the relation given by the degree of a polynomial. Because of this, it is usual to call congruences as \textbf{graduations} when we deal with the semigroup of polynomials or rational functions.\\
Given two congruences $\mathcal{A}=\{A_i\}_{i\in I}$ and $\mathcal{B}=\{B_j\}_j\in J$, we say that $\mathcal{A}$ and $\mathcal{B}$ are equivalents if 
$$A_i\cap B_j\neq\emptyset, \text{ for some }i\in I,\ j\in J,\qquad\Rightarrow\qquad A_i=B_j.$$
In general, given two congruences $\mathcal{A}=\{A_i\}_{i\in I}$ and $\mathcal{B}=\{B_j\}_j\in J$, we say that congruence $\mathcal{A}$ is subordinated to congruence $\mathcal{B}$, which we denote by $\mathcal{A}\leq\mathcal{B}$, if for every $i\in I$ there exists $j\in J$ such that 
$$A_i\subset B_j.$$
Note that two congruences are equivalents if and only if they are subordinated to each other. \\
If $\varphi:S\rightarrow \CC$ is a multiplicative character over a semigroup S, then \textit{the fibers of the character} define a graduation over $S$, denoted by $\mathcal{P}_{\varphi}$. \\
We are in position to give another equivalence for an elementary character. 
\begin{prop}
Let $\varphi:S\rightarrow \CC$ be a character over $S$, with $S\in\{\pol,\rat\}$. Suppose that $\mathcal{P}_{\varphi}$ is subordinated to $\mathcal{P}_0$, the canonical graduation. Then $\varphi$ is elementary, that is, there  exists a multiplicative complex-value function $f:\NN\rightarrow\CC$ such that 
$$\varphi(R)=f\circ\emph{deg}(R).$$ 
\end{prop}
\bproof Given $\varphi:S\rightarrow\CC$, define the function $f:\NN\rightarrow\CC$ as
$$f(n)=\varphi(x^n).$$
Note that $f$ is well defined and 
$$f(n\cdot m)=\varphi(x^{nm})=\varphi((x^m)^n)=\varphi(x^n\circ x^m)=\varphi(x^n)\varphi(x^m)=f(n)f(m).$$ 
Now, since $\mathcal{P}_{\varphi}$ is subordinated to $\mathcal{P}_0$ every element $s\in S$ is related to the polynomial $x^{\emph{deg}(s)}$, then 
$$\varphi(R)=\varphi(x^{\emph{deg}(R)})=f(\emph{deg}(R)),$$
so $\varphi$ is elementary.  \eproof

This way, we say a character $\varphi:S\rightarrow\CC$ where $S\in\{\pol,\rat\}$ is elementary if the graduation $\mathcal{P}_{\varphi}$ is subordinated to the canonical graduation. \\

In a natural way, the set of multiplicative (resp. additive) characters over a semigroup $S$, form a semigroup under multiplication (resp. addition), that is, if $\varphi$ and $\psi$ are multiplicative characters over a semigroup $S$, it is possible define a new multiplicative character $\chi$over $S$ as 
$$\chi(s)=\varphi(s)\cdot\psi(s),\ s\in S.$$ 
The case of additive characters is defined analogously as the sum. If $\chi$ is defined as the product of $\varphi$ and $\psi$, is there any relation between $\mathcal{P}_{\chi}$, $\mathcal{P}_{\varphi}$ and $\mathcal{P}_{\psi}$? We recall that the intersection of congruences is again a congruence. 
\begin{prop}
Given multiplicative characters $\varphi$ and $\psi$ over a semigroup $S$, if $\chi$ is the product character of $\varphi$ and $\psi$, $\chi=\varphi\cdot\psi$, then $\mathcal{P}_{\varphi}\cap \mathcal{P}_{\psi}\leq \mathcal{P}_{\chi}$. 
\end{prop}
\bproof Denote $\mathcal{P}_I=\mathcal{P}_{\varphi}\cap\mathcal{P}_{\psi}$ and suppose $s\mathcal{P}_It$, that is, $\varphi(s)=\varphi(t)$ and $\psi(s)=\psi(t)$, then
$$\chi(s)=\varphi(s)\cdot\psi(s)=\varphi(t)\cdot\psi(t)=\chi(t),$$
hence $s\mathcal{P}_{\chi}t$. \eproof

In some cases, depending on the nature of the characters, it is possible to have equivalence among the characters instead of only subordination. Take as an example the semigroup $S=\pol$ of complex polynomials and the characters $\varphi=\emph{deg}$ and $\psi=e^l$, where $l=$\textit{length of the prime decomposition}. So we define the product as
$$\chi(P)=\text{deg}(P)\cdot e^{l(P)},\ \text{for } P\in\pol.$$ 
Now, suppose $\chi(P)=\chi(Q)$ for some $P,Q\in\pol$, then 
$$\chi(P)=\text{deg}(P)\cdot e^{l(P)}=\text{deg}(Q)\cdot e^{l(Q)}=\chi(Q),$$
hence
$$e^{l(Q)-l(P)}=\frac{\text{deg}(P)}{\text{deg}(Q)}\in\QQ,$$
since $l(Q),l(P)\in\NN$, then $e^{l(Q)-l(P)}\in\QQ\ \Leftrightarrow\ l(Q)=l(P)$, which implies that $\text{deg}(P)=\text{deg}(Q)$, so $\mathcal{P}_{\chi}\leq\mathcal{P}_{\varphi}\cap\mathcal{P}_{\psi}$, and by above proposition $\mathcal{P}_{\varphi}\cap\mathcal{P}_{\psi}\leq\mathcal{P}_{\chi}$, then both graduations are equivalent.
\section{Examples given rise to conjecture}
In the present section, we present some examples concerning \textbf{Conjecture A} and general setting of possible obstructions to the \textit{polynomial character extension problem}. 
\subsection{Restrictions over $Aff(\CC)$}
Suppose $\varphi$ is a multiplicative character on complex polynomials such that $\varphi\large|_{Aff(\CC)}$ is non constant and $\varphi(\langle\alpha z\rangle)$ is not identically $1$, if $\alpha$ is and $n-$root of unit, then
$$\varphi(z^n)=\varphi(z^n\circ \alpha z)=\varphi(z^n)\cdot\varphi(\alpha z)\Rightarrow \varphi(z^n)=0.$$
Similar calculations can be made for other polynomials $P$ with $\text{deg}(P)\geq2$. Hence, we define the followings restrictions over the group $Aff(\CC)$ for a given character $\varphi$ on polynomials.
\begin{rest}\label{rst}
Given a character $\varphi:\pol\rightarrow \CC$, the following conditions are considered.
\begin{itemize}
\item[1.] $\varphi(c)\equiv0$ for every constant $c\in\CC$ and
\item[2.] $\varphi(A(z))=1$ for every $A\in Aff(\CC)$.  
\end{itemize}
\end{rest}
Since $PSL(2,\CC)=\langle Aff(\CC),I(z)\rangle$, where $I(z)=z^{-1}$, the above restrictions extend to $PSL(2,\CC)$ with 
\begin{itemize}
\item[3.] $\varphi(z^{-1})=\pm1$, since $1=\varphi(z)=\varphi(1/z\circ 1/z)=\varphi(1/z)\varphi(1/z)$. 
\end{itemize}
\subsection{Example zero}
As was mention above, in \cite{cmp}, authors constructed a non trivial character over polynomials, which is non elementary. This character will serves as the first brick in the construction of examples for Conjecture A. \\
Recall that by the first part or Ritt's Theorem, for a given polynomial $P\in\pol$, the function given by 
\begin{displaymath}
\begin{array}{rcl}
\mathit{L}:\pol & \rightarrow & \CC\\
			P & \mapsto & \mathit{L}(P)=e^{l(P)},
\end{array}
\end{displaymath}
where $l(P)$ is the length of prime decompositions of $P\in\pol$, defines a multiplicative character. Since $l(P)$ do not depends on the grade of $P$, $\mathit{L}$ is non elementary. \\
The following simple example as given in \cite{bere} and \cite{ber}, shows that this character cannot be extended to rational functions. \\
Consider the function 
$$R(z):=\left(\frac{z^4-8z}{z^3+2z^2+2z+1}\right)^3.$$
Bergweieler proves that this function has the following two prime decompositions:
\begin{equation}\label{ex2}
R(z)=z^3\circ \frac{z^2-4}{z-1}\circ\frac{z^2+2}{z+1}=\frac{z(z-8)^3}{(z+1)^3}\circ z^3.
\end{equation}
This double decompositions makes the function $l(R)$ not well defined, and then the multiplicative character $\mathit{l}$ cannot be extended to rational functions. \\
This example is sometimes called as the \emph{Ritt's example} of a non extendable character. 
\subsection{Example in $z^p$} Consider the following preliminaries for the next example. 
First of all, for a given rational map $R$, consider its \textit{Hurwitz class} $H(R)$, defined as 
$$H(R)=\{Q\in\rat|\exists\phi,\psi\in\text{Homeo}(\overline{\CC}) \text{ with }Q\circ\phi=\psi\circ R\},$$
that is, the set of all rational maps with the same combinatorics for the first iteration. If $\phi=\phi(R)$ and $\psi=\psi(R)$, in the above definition, are even \textit{conformal}, then the space is called the \textit{conformal Hurwitz class} and denoted by $CH(R)$. Denote by $Crit(P)$ the set of critical points of a polynomial map $P$. Since $|Crit(z^n)=2|$, $n\geq3$ (where are considering the point at infinity as a point of the Riemann sphere $\overline{\CC}$), and $|Crit(T_n)|=n$ where $T_n$ is the general \textit{Tchebycheff polynomial}, then it is clear that 
$$CH(z^n)\neq CH(T_n).$$
Finally, we recall the second part of the of the Ritt's Theorem \cite{ritt1}. 
\begin{theorem}\label{decom}
Any pair of prime decompositions of a polynomial $P$ are related by a finite number of Ritt's transformations. Moreover, there are three types of Ritt's transformations
\begin{itemize}
\item[1.] $P_i\circ P_{i+1}\leftrightarrow (P_i\circ A)\circ (A^{-1}\circ P_{i+1})$ where $A\in Aff(\CC)$.
\item[2.] $P_i\circ P_{i+1}\leftrightarrow P_{i+1}\circ P_i$ where $P_j$ are Tchebychev polynomials. 
\item[3.] $z^k\circ z^rQ(z^r)\leftrightarrow z^k(Q(z))^k\circ z^r$, here $Q(z)\in\pol$. 
\end{itemize}
\end{theorem}
We are able to define the following character over $\pol$. For a fixed prime number $p\geq3$,
\begin{displaymath}
\varphi_1(P) = \left\{ \begin{array}{ll}
p & \textrm{if $P\in H(z^p)$}\\
1 & \textrm{otherwise.}
\end{array} \right. 
\end{displaymath}
For the above reasoning, the Theorem \ref{defcha}, and the second part of the Ritt's Theorem (Theorem \ref{decom}), it is clear that $\varphi_1$ is well defined and it is a \textit{multiplicative character} over $\pol$. Moreover, $\varphi_1$ is not elementary. Recall that $\varphi_1$ holds Character Restrictions \ref{rst}. 

If we try to extend $\varphi_1$ to rational functions $\rat$, we most to consider the fact that every rational function $f(z)$ of degree 2 is in the Hurwitz class of $z^2$. Also, it is well known that for every Tchebychev polynomial $T_p$ of degree $p\geq3$ the following relation holds:
$$T_p\circ Y(z) = Y\circ z^p, \text{ where }Y(z)=\frac{1}{2}(z+z^{-1}).$$
Applying $\varphi_1$ to the above relation, we conclude that $\varphi_1(Y(z))=\varphi_1(z^2)=0$ which contradicts the definition of $\varphi_1$. Then $\varphi_1$ cannot be extended to rational functions $\rat$. 
\subsection{Example in grade 4}
One of the significant difference between the semigroup of complex polynomials $\pol$ and the semigroup of complex rational functions is the Ritt's Theorem. Several authors, see for example \cite{ber} and \cite{pak}, have exhibit examples of rational functions with prime decomposition of different length. Klein functions, see \cite{kl}, provide the simplest examples of rational functions for which Ritt's Theorem fails to be true. 

First, we define the following character over $\pol$. 
\begin{displaymath}
\varphi_2 = \left\{ \begin{array}{ll}
4 & \textrm{if $P\in \mathcal{P}ol_4(\CC)\backslash H(z^4)$}\\
1 & \textrm{otherwise.}
\end{array} \right.
\end{displaymath}
Again, with the same reasoning in the above section and the second part of Ritt's Theorem, $\varphi_2$ is well defined and is a \textit{multiplicative character} over $\pol$. To see that $\varphi_2$ cannot be extended to rational functions, consider the following function and its decompositions
$$-\frac{1}{432}\frac{(16z^8-56z^4+1)^3}{z^4(4z^4+1)4}=\left(\frac{1}{54}\frac{(z+7)^3}{(z-1)^2}\right)\circ\left(\frac{1}{2}\left(z+\frac{1}{z}\right)\right)\circ(-z^2)\circ z^2$$
and 
$$-\frac{1}{432}\frac{(16z^8-56z^4+1)^3}{z^4(4z^4+1)4}=\left(-\frac{256}{27}z^3(z-1)\right)\circ\left(\frac{1}{4}\frac{(z-1)^3}{z^2+1}+1\right)\circ\left(z-\frac{1}{2z}\right).$$
These decompositions, correspond to the chains 
$$1<C_2<C_4<D_8<S_4,\qquad 1<C_2<S_3<S_4$$
of the group $S_4$, which is the \textit{monodromy group} (see Appendix in \cite{pak}) of the function 
\begin{equation}\label{ex4}
f_{S_4}(z)=-\frac{1}{432}\frac{(16z^8-56z^4+1)^3}{z^4(4z^4+1)4}.
\end{equation}
Note that, for the functions of degree 3, we have the following chains of transformations. (\emph{Maps over the arrows are post-composed and maps under the arrows are pre-composed.}) 

$$\left(\frac{1}{4}\frac{(z-1)^3}{z^2+1}+1\right)\xrightarrow[]{4(z-1)}\frac{(z-1)^3}{z^2+1}\xrightarrow[z+1]{}\frac{z^3}{(z+1)^2+1}\xrightarrow[]{1/z}$$
$$\xrightarrow[]{1/z}\frac{(z+1)^2+1}{z^3}\xrightarrow[1/z]{}\frac{(1/z+1)^2+1}{(1/z)^3}=\frac{z^3[(1+z)^2+z^2]}{z^2}=z[(1+z)^2+z^2],$$
and
$$\left(\frac{1}{54}\frac{(z+7)^3}{(z-1)^2}\right)\xrightarrow[]{54z}\left(\frac{(z+7)^3}{(z-1)^2}\right)\xrightarrow[-z]{}\frac{(7-z)^3}{(z+1)^2}\xrightarrow[]{1/z}\frac{(z+1)^2}{(7-z)^3}\xrightarrow[z+7]{}$$
$$\xrightarrow[z+7]{}\frac{(z+8)^2}{-z^3}\xrightarrow[1/z]{}\frac{(1/z+8)^2}{-1/z^3}=\frac{-z^3(1+8z)^2}{z^2}=-z(1+8z)^2.$$
This way, we have 
$$\left(\frac{1}{4}\frac{(z-1)^3}{z^2+1}+1\right),\left(\frac{1}{54}\frac{(z+7)^3}{(z-1)^2}\right)\in\pol, $$
which implies that 
$$\varphi_2\left(\frac{1}{4}\frac{(z-1)^3}{z^2+1}+1\right)=\varphi_2\left(\frac{1}{54}\frac{(z+7)^3}{(z-1)^2}\right)=1.$$
Hence 
$$\varphi_2\left(\left(\frac{1}{54}\frac{(z+7)^3}{(z-1)^2}\right)\circ\left(\frac{1}{2}\left(z+\frac{1}{z}\right)\right)\circ(-z^2)\circ z^2\right)=1$$
but
$$\varphi_2\left(\left(-\frac{256}{27}z^3(z-1)\right)\circ\left(\frac{1}{4}\frac{(z-1)^3}{z^2+1}+1\right)\circ\left(z-\frac{1}{2z}\right)\right)=4,$$
then $\varphi_2$ is not well defined over $\rat$, so it cannot be extended.

We had constructed two characters, well defined over polynomials, that cannot be extended to rational functions. 
\section{Ideal extensions and some notes on obstructions}
Given semigroups $S$, $E$, and a semigroup with zero $Q$, we say that $E$ is an \textbf{ideal extension} of $S$ by $Q$ if $S\subset E$, $S$ is an ideal of $E$ and $E/S=Q$. Trivially, every semigroup $S$ is an ideal extension of every proper ideal $I\subset S$, with Rees quotient $Q=S/I$. 

Since the second example of previous sections, the fact that, in general, rational functions do not holds Ritt's Theorem seems to be an obstructions to the extension of characters over $\pol$. Note that the set of rational functions for which Ritt's Theorems fails to be true, form an ideal of $\rat$:
$$I_R=\{Q\in\rat| Q \text{ has prime decompositions with different length}\}.$$
We call this, the \textbf{Ritt ideal} of $\rat$. Consider first the ideal of decomposable rational functions
$$I_D=\{Q\in\rat:Q\text{ is decomposable}\}.$$
In this case, the Rees quotient for $I_D$ is given by the set of prime rational functions. Moreover, every element of this semigroup with zero $Q_D=\rat/I_D$ is an idempotent element
$$Q\in\rat\backslash I_D\Rightarrow Q\circ Q\in I_D.$$
Suppose there is a character defined over $I_D$, $\psi:I_D\rightarrow \CC$: Can $\psi$ be defined for every $R\in Q_D^*$, extending the character $\psi$?. 

Given a fix $R\in Q^*$, $R^n\in Q^*$ for  every $n\geq2$. Then, $\psi$ is well defined over $\{R^n\}_{n\geq2}$. Moreover, since $R^6=R^2\circ R^2\circ R^2=R^3\circ R^3$ we have $\psi(R^2)^3=\psi(R^3)^2$. So, we define $\psi(R)=\psi(R^3)/\psi(R^2)$. Then $\psi$ is well define over the semigroup $\{R^n\}_{n\geq1}$ and is a character. 

\begin{lemma}\label{itera}
Let $I\subset \rat$ be an ideal and $\varphi:\rat\rightarrow\CC$ a complex value function such that $\varphi\large|_{I}$ is a multiplicative character. If $R\notin I$ but $R^2\in I$, then $\varphi$ is a character over $\langle I,R\rangle$ with
$$\varphi(R)=\frac{\varphi(R^3)}{\varphi(R^2)}.$$
\end{lemma}
\bproof Since $I$ is an ideal is enough to prove that $\varphi(R\circ Q)=\varphi(R)\varphi(Q)$ and $\varphi(Q\circ R)=\varphi(Q)\varphi(R)$ for every $Q\in I$. Note that $\{R^n\}_{n\geq2}\subset\rat$ ($I$ is ideal) and then $\varphi(R)=\varphi(R^3)/\varphi(R^2)$ is well defined.  Given that $\varphi$ is a character over $I$, we have
$$\varphi(R^3)\varphi(Q)=\varphi(R^3\circ Q)=\varphi(R^2\circ R\circ Q)=\varphi(R^2)\varphi(R\circ Q),$$
from which we obtain $$\varphi(R\circ Q)=\frac{\varphi(R^3)}{\varphi(R^2)}\varphi(Q),$$
and then the first equality holds. For the second equality the proof is analogous.\eproof
\begin{theorem}
Given a complex valued function $\varphi:\rat\rightarrow\CC$ such that $\varphi\large|_{I_D}$ is a multiplicative character, then $\varphi:\rat\rightarrow\CC$ can be extended as a multiplicative character.
\end{theorem}
\bproof Note that $I_D$ holds hypothesis in the above lemma. Take $R_0\notin I_D$, then $\varphi$ is a character over $\langle I_D,R_0\rangle$. Also, $I_{D0}=\langle I_D,R_0\rangle\leq\rat$ is an ideal with $I_{D0}\backslash I_{D}=R_0$ satisfying hypothesis in Lemma \ref{itera}. Then for $R_1\notin I_{D0}$, $\varphi$ is a character over $I_{D1}=\langle I_{D0},R_1\rangle$ which again is an ideal with $I_{D1}\backslash I_{D0}=R_1$. This argument can be repeated indefinitely. Although $\rat$ is not isomorphic to $\mathbb{N}$, and we cannot use an induction argument, we can fulfilled $\rat$ via $I_D$ in this way and then use a transfinite limit argument to extend $\varphi$ over all of $\rat$, and then $\varphi$ is a character over $\rat$. \eproof

In general terms, we can proved the following on extension over ideal of semigroups.
\begin{theorem}\label{idealextension}
Let $\varphi:S\rightarrow\CC$ be a well defined complex function over a semigroup $S$ and let $I\subset S$ be an ideal. If $\varphi\large|_I$ is a multiplicative character, then, $\varphi$ can be extended to a character $\widetilde{\varphi}$ over $S$. 
\end{theorem}
\bproof First, we construct the possible extension. Let $R\in S\backslash I$, for every $Q\in I$ with $\varphi(Q)\neq0$, we define
$$\phi_Q(R):=\frac{\varphi(R\cdot Q)}{\varphi(Q)}.$$
\textbf{CLAIM 1}. For $R\in S\backslash I$ and $Q\in I$, 
$$\varphi(Q\cdot R)=\varphi(R\cdot Q),$$
and then $\phi_Q$ is well defined. \\
Note that $Q\cdot R\cdot Q,\ R\cdot Q,\ Q\cdot R\in I$ since $I$ is an ideal, and by hypothesis $\varphi\large|_I$ is a character, if $\varphi(Q)\neq0$, then we have
$$\varphi(R\cdot Q)=\frac{\varphi(Q)\cdot\varphi(R\cdot Q)}{\varphi(Q)}=\frac{\varphi(Q\cdot R\cdot Q)}{\varphi(Q)}=\frac{\varphi(Q\cdot R)\cdot\varphi(Q)}{\varphi(Q)}=\varphi(Q\cdot R).$$
If $\varphi(Q)=0$, then we define $\varphi(R\cdot Q)=0=\varphi(Q\cdot R)$. \\
\textbf{CLAIM 2}. For each $R\in S\backslash I$, $\phi_Q$ is constant. \\
Let $Q_1,\ Q_2\in\ I$, be such that $\varphi(Q_1)\neq0\neq\varphi(Q_2)$, then 
\begin{displaymath}
\begin{array}{rcl}
\phi_{Q_1}(R)-\phi_{Q_2}(R) & = & \frac{\varphi(Q_1\cdot R)}{\varphi(Q_1)}-\frac{\varphi(Q_2\cdot R)}{\varphi(Q_2)} \\ \\
						& = & \frac{\varphi(Q_1\cdot R)\varphi(Q_2)-\varphi(Q_2\cdot R)\varphi(Q_1)}{\varphi(Q_1)\varphi(Q_2)} \\ \\
						& = & \frac{\varphi(Q_1\cdot R)\varphi(Q_2)-\varphi(Q_2\cdot R\cdot Q_1)}{\varphi(Q_1)\varphi(Q_2)} \\ \\
						& = & \frac{\varphi(Q_1\cdot R)\varphi(Q_2)-\varphi(Q_2)\varphi(R\cdot Q_1)}{\varphi(Q_1)\varphi(Q_2)} \\ \\
						& = & \frac{\varphi(Q_1\cdot R)-\varphi(R\cdot Q_1)}{\varphi(Q_1)}  \\
						& = & 0 \textit{  ( by claim 1).}
\end{array}
\end{displaymath}
We define the extension of $\varphi$ as 
\begin{displaymath}
\widetilde{\varphi}(R) = \left\{ \begin{array}{ll}
\varphi(R) & \textrm{if $R\in I$}\\
\phi_Q & \textrm{otherwise, for some $Q$ with $\varphi(Q)\neq0$}.
\end{array} \right.
\end{displaymath}
We have to prove that $\widetilde{\varphi}(R\cdot Q)=\widetilde{\varphi}(R)\cdot\widetilde{\varphi}(Q)$, for every $R,Q\in S$, we have to consider three cases.
\begin{itemize}
\item $R,Q\in I$: in this case $\widetilde{\varphi}(R\cdot Q)=\varphi(R\cdot Q)=\varphi(R)\cdot \varphi(Q)=\widetilde{\varphi}(R)\cdot\widetilde{\varphi}(Q)$. 
\item $R\in S\backslash I$ and $Q\in I$: then $R\cdot Q\in I$, so 
$$\widetilde{\varphi}(R\cdot Q)=\varphi(R\cdot Q)=\frac{\varphi(R\cdot Q)\cdot\varphi(Q)}{\varphi(Q)}=\widetilde{\varphi}(R)\cdot\widetilde{\varphi}(Q).$$
\item $R,Q\in S\backslash I$ and $R\cdot Q\in S\backslash I$: for this, given $S_1\in I$, with $\varphi(S_1)\neq0$, consider the following
\begin{displaymath}
\begin{array}{rcl}
\widetilde{\varphi}(R\cdot Q) & = & \frac{\varphi(R\cdot Q\cdot S_1)}{\varphi(S_1)}\\ \\
							  & = & \frac{\widetilde{\varphi}(R)\cdot \varphi(Q\cdot S_1)}{\varphi(S_1)}\\ \\
							  & = & \frac{\widetilde{\varphi}(R)\cdot\widetilde{\varphi}(Q)\cdot\varphi(S_1)}{\varphi(S_1)}\\ \\
							  & = & {\widetilde{\varphi}(R)\cdot\widetilde{\varphi}(Q)}
\end{array}
\end{displaymath}

\end{itemize}
This proves the proposition. \eproof \\

\begin{cor}
There exists a 1-1 correspondence between characters over $\rat$ and proper ideals of $\rat$ with a well defined character.
\end{cor}
\bproof Since Theorem \ref{idealextension}, if $I\subset\rat$ is an ideal with $\varphi:I\rightarrow \CC$ a character then $\varphi$ is extended to a character $\widetilde{\varphi}:\rat\rightarrow\CC$. \\
On the other hand, let $\psi:\rat\rightarrow\CC$ be a character and $\varphi:I\rightarrow$ be a character over a proper ideal $I\subset\rat$ with $\psi\large|_I=\varphi$, we need to prove that $\psi=\widetilde{\varphi}$.\\
Let $R\in\rat\backslash I$ and $Q\in I$ arbitrary with $\varphi(Q_0)\neq0$, by definition 
$$\widetilde{\varphi}(R)=\frac{\varphi(R\circ Q)}{\varphi(Q)},$$
now, since $\psi$ is a character, we have $\psi(R\circ Q)=\psi(R)\psi(Q)$, hence
$$\psi(R)=\frac{\psi(R\circ Q)}{\psi(Q)}=\frac{\varphi(R\circ Q)}{\varphi(Q)}=\widetilde{\varphi}(R).$$
Then $\psi=\widetilde{\varphi}$. \eproof \\

Consider now the case $I=I_R$, the Ritt ideal. In this case $Q_{I_R}=\rat/I_R$. Proposition \ref{idealextension}, tell us that no matter how the ideal is, the extension is possible. Then, if we have a character defined over $I_R$, this can be extended to $\rat$. Some questions arise with this: Is is possible to define a non elementary character over $I_D$? And, if this character exists, its restriction to $\pol$ is again non elementary?
\subsection{A cyclic character}
Under Character Conditions \ref{rst}, we study a special character defined as follows. Consider the set $$\textbf{P}=\{R\in\rat:R\text{ is not decompsosable}\},$$ then we define a function 
$$\delta:\textbf{P}\rightarrow\RR_+$$ 
such that $\delta\large|_{\textbf{P}}\in\{\alpha,\alpha^k\}$ for some $\alpha\in\RR_+\backslash\{0,1\}$ and $2\leq k\in\NN$, i.e., if $R\in\textbf{P}$, then $\delta(R)=\alpha$ or $\delta(R)=\alpha^k$. Moreover, suppose $\delta\large|_{\textbf{P}\cap\mathcal{R}at_d}$ is constant for every $d\geq2$, where $\mathcal{R}at_d$ is the set of rational functions of degree $d$. We would like to see if $\delta$ may defined a character over $I_R$, the ideal of rational function with prime decompositions of different length. \\
Let $R\in I_R$, with prime decompositions
$$R=R_1\circ R_2\circ...\circ R_i=S_1\circ S_2\circ...\circ S_j,\ (1<i<j).$$
To $\delta$ be a character, we must have 
\begin{displaymath}
\begin{array}{rcccl}
\delta(R) & = & \delta(R_1\circ R_2\circ...\circ R_i) & = & \delta(S_1\circ S_2\circ...\circ S_j) \\
          & = & \delta(R_1)\cdot\delta(R_2)\cdot...\cdot\delta(R_i) & = & \delta(S_1)\cdot\delta(S_2)\cdot...\cdot\delta(S_j) 
\end{array}
\end{displaymath} 
implying equations
\begin{displaymath}
\begin{array}{rcl}
\alpha^{r_1}\cdot\alpha^{r_2}\cdot...\cdot\alpha^{r_i} & = & \alpha^{s_1}\cdot\alpha^{s_2}\cdot...\cdot\alpha^{s_j} \\
\alpha^{r_1+r_2+...+r_i} & = & \alpha^{s_1+s_2+...+s_j} \\
 r_1+r_2+...+r_i & = & s_1+s_2+...+s_j
\end{array}
\end{displaymath}
where $r_{i'},s_{j'}\in\{1,k\}$, $1\leq i'\leq i$ and $1\leq j'\leq j$.  Note that if $r_{i'}=s_{j'}=1$ for every $i'$ and $j'$, then $i=j$ which contradicts the assumption $R\in I_R$. Analogously if $r_{i'}=s_{j'}=k$. Suppose then that there exist different $r_{i'}'s$ and/or different $s_{j'}$. \\
Let $p=|\{r_{i'}=1\}|$, $q=|\{r_{i'}=k\}|$ and $m=|\{s_{j'}=1\}|$, $n=|\{s_{j'}=k\}|$, we have the following relations: 
$$p+q=i,\qquad m+n=j,\qquad\text{and}\qquad p+qk=m+nk.$$
Since $k\geq2$, 
$$q+qp=p+q+q(k-1)=i+q(k-1)=m+nk=m+n+n(k-1)=j+n(k-1),$$
then
$$(k-1)=\frac{j-i}{n-q},$$
but $k\in\NN$, this implies that $(k-1)\in\NN$ so $\frac{j-i}{n-q}\in\NN$ with $k=1+\frac{j-i}{n-q}$. Given \textit{Example Zero} (\ref{ex2}) above, $j-i=3-2=1$ and the fact that $k$ is constant, $k=2$. So, in the general case, we have the equalities
$$j-i=q-n \qquad\text{and}\qquad m-p=2(j-i).$$
We obtained the following restrictions:
\begin{itemize}
\item If $q=0$ then $n=0$ which implies $i=j$, a contradiction, so $q>0$;
\item if $m=0$ then $p=0$ which implies $i=j$, again a contradiction so $m>0$;
\item if $n=0$ then $q=i-j$ and if $p=0$ then $m=2(j-i)$. 
\end{itemize}
Moreover, $\delta$ may not be \textit{elementary}: Suppose $R=R_1\circ R_2\circ R_3$, with $\deg R=d=d_1d_2d_3$, $\deg R_i=d_i$. Then $\log_{\alpha}\delta(R)\geq3$, now let $Q\in\rat$ with $\deg Q=d$ and $Q\in\textbf{P}$, then $\log_{\alpha}(Q)\leq2$, so $\deg(R)=d=\deg (Q)$ but $\delta(R)\neq\delta(Q)$ .  \\
On the other hand, applying $\delta$ to rational function (\ref{ex2}), we have that $\delta\large|_{\textbf{P}\cap\mathcal{R}at_4}=\alpha^2$. Suppose there exists a rational function $R$ of degree $d=24$ with decompositions
$$R=R_1\circ R_2\circ R_3=S_1\circ S_2,$$
with degrees $\deg S_1=8$ and $\deg S_2=3$ (probably not in that order), and $\deg R_i=4$ for some $i=1,2,3$. Then $\delta(S_1)=\alpha^{3}$, which contradicts the definition of $\delta$. Unfortunately, there is no known example of this type, and in general there are only few examples of functions with prime decompositions of different length that may yield this contradiction. In other words, only few \textit{generators} elements of ideal $I_R$ are known, two of them are functions (\ref{ex2}) and (\ref{ex4}). \\
Apart from $\pol$, the so called \textit{Laurent Polynomials} $\mathcal{L}\rat$ also holds Ritt's Theorem. A question arise: Is $Q_{I_R}=\pol\sqcup\mathcal{L}\pol$? That is, are the polynomials and Laurent polynomials, the only subsemigroups of $\rat$ that holds Ritt's Theorem? 
\section{Addendum}
Let $B$ be a rational function. As an element of the semigroup of (complex) rational functions, F. Pakovich in \cite{Pa2}, consider a \textit{graph}, based on $B$, which gives a topological realization of the decomposition structure of $B$. \\
Given a decomposition $B=V\circ U$ of a rational function $B$, where $U,V$ are non constant rational functions, the function $\hat{B}=U\circ V$ is called an \textit{elementary transformation} of $B$. We consider the following relation: we say that $A\sim B$ if there exists a chain of elementary transformations from $B$ to $A$, in other words, there exist rational functions $U_i$, $V_i$, $1\leq i\leq s$, and the chain
$$B\mapsto B_1\mapsto B_2\mapsto \dots \mapsto B_s=A,$$
where 
\begin{center}
$B=V_1\circ U_1,\ B_1=U_1\circ V_1,\ B_s=U_s\circ V_s=A, and$ \\ $B_i=U_i\circ V_i=V_{i+1}\circ U_{i+1},\ 2\leq i\leq s-1.$
\end{center}
Note that in the above definition, there is no restriction on the degree of maps $U$ and $V$, so the equivalence class [B] of a rational functions is a union of conjugacy classes (we may write $B=(B\circ\mu^{-1})\circ\mu$, where $\mu$ is a M\"obius transformation). This way, it is possible to the define a multigraph $\Gamma_B$, with vertices given by fixed representatives of conjugacy classes in $[B]$, and edges, connecting $B_i$ with $B_j$, given by solutions of the system
\begin{equation}\label{edges}
B_i=V\circ U,\qquad B_j=U\circ V, 
\end{equation}
in rational functions. Note that loops in $\Gamma_B$ may exists, they correspond to solutions of 
\begin{equation}\label{loops}
B_i=V\circ U=U\circ V. 
\end{equation}
Some immediate results about $\Gamma_B$ are listed below (see \cite{Pa2}).
\begin{lemma}
The graph $\Gamma_B$ does not depends on the choice of representatives of conjugacy classes in $[B]$. 
\end{lemma}
\begin{theorem}
Let $B$ be a rational function of degree at least two. Then the graph $\Gamma_B$ is finite, unless $B$ is \emph{flexible Latt\`es map}. 
\end{theorem}
Few basic but interesting examples are exhibit in \cite{Pa2}. \\
Given the nature of this work, we are interested in results concerning rational functions with prime decompositions of different length. In this direction, F. Pakovich consider the following result.\\
Recall that, two decompositions (maximal or not) of a given rational functions $R$, having the same length 
$$R=R_r\circ R_{r-1}\circ...\circ R_1,\qquad R=S_r\circ S_{s-1}\circ...\circ S_1,$$
are called equivalent if either $r=1$ and $R_1=S_1$ or $r\geq2$ and there exist M\"obius transformations $\mu_i$, $1\leq r\leq r-1$, such that
$$R_r=S_r\circ\mu_{r-1},\qquad R_i=\mu_i^{-1}\circ S\circ \mu_{i-1},\qquad 1<i<r, \text{ and }R_1=\mu_1^{-1}\circ S_1.$$
Then, given a rational function $B$ with a maximal decomposition 
\begin{equation}
B=U_r\circ U_{r-1}\circ...\circ U_1,
\end{equation}
we say that $B$ is \textit{generically decomposable} if 
\begin{itemize}
\item each function 
$$\widehat{B}_i=(U_i\circ U_{i-1}\circ...\circ U_1)\circ(U_r\circ U_{r-1}\circ...\circ U_{i+1}),$$
has a unique equivalence class of maximal decompositions, and
\item the functions $\widehat{B}_i$, $1\leq i\leq r-1$ are pairwise not conjugated.
\end{itemize}
Finally, we define the graph $\Gamma^0_B$ as the subgraph of $\Gamma_B$ with all the loops, corresponding to automorphisms, removed. We say that the rational function $B$ is not \textit{special} if is neither a Latt\`es map, nor conjugated to $z^{\pm}$ or $T_n$, with $T_n$ a Tchevychev polynomial. 
\begin{theorem}[Lemma 6.2, \cite{Pa2}]
Assume that a non-special rational function $B$ having a maximal decomposition of length $r$ is generically decomposable. Then $\Gamma_B^0$ is the complete graph $K_r$. 
\end{theorem}
In terms of the present work, for a rational function $B$, be generically decomposable means that the function $B$ and all of the ``\textit{shifted}" functions (of its factors) $\widehat{B}_i$ hold First Ritt's Theorem, in other words, they do not belong to the Ritt Ideal $I_R$ defined above. 
As was noticed, the Ritt Ideal is an obstruction for the extension of the \textit{length character}, non elementary, well define for the semigroup of polynomials. It would be important to have a characteristic to distinguish rational functions belonging to the Ritt Ideal. In this direction, the following conjecture is proposed. 
\begin{conje}
Given a non-special rational function $B$, if the graph $\Gamma_B\backslash{B}$ is disconnected, then $B$ belongs to the Ritt Ideal. 
\end{conje}

\end{document}